\documentclass[12pt]{article}
\usepackage[utf8]{inputenc}
\usepackage{amsmath}
\usepackage{amsfonts}
\usepackage{amssymb}
\usepackage{amsthm}
\usepackage{relsize}
\usepackage{array}
\usepackage{multirow}
\usepackage{tabu}
\usepackage{soul}
\usepackage{graphicx}
\usepackage{epigraph}
\usepackage{float}
\usepackage[english]{babel}
\usepackage[usenames, dvipsnames]{color}
\usepackage{fancyhdr}
\usepackage[Sonny]{fncychap}
\usepackage{tikz-cd}
\usepackage{tkz-euclide}
\usepackage{sseq}
\usepackage{newfloat}
\usepackage{amsmath,mathtools}
\usepackage{pgfplots}
\usepackage{geometry}
 \usepackage{sectsty}
\sectionfont{\fontsize{14}{15}\selectfont}

\usepackage[all]{xy}
\DeclareFloatingEnvironment[fileext=dia]{diagram}

\theoremstyle{definition}

\usepackage[most]{tcolorbox}

\makeatletter
\newtcolorbox{note}[1][]{%
	breakable,
	enhanced jigsaw, 
	borderline west={3pt}{0pt}{black!10!white}, 
	borderline south={1pt}{0pt}{black!10!white}, 
	borderline east={1pt}{0pt}{black!10!white},
	borderline north={1pt}{0pt}{black!10!white},
	sharp corners, 
	boxrule=0pt, 
	attach title to upper, 
	left=0pt,
	right=0pt,
	top=0pt,
	bottom=0pt,
	boxsep=5pt,
	colback=white,
	frame hidden,
	#1
}
\makeatother

\makeatletter
\newtcolorbox{note1}[1][]{%
	breakable,
	enhanced jigsaw, 
	sharp corners, 
	boxrule=0pt, 
	attach title to upper, 
	fontupper=\linespread{1.1}\fontfamily{qpl}\selectfont,
	fontlower=\linespread{1.1}\fontfamily{qpl}\selectfont, 
	left=0pt,
	right=0pt,
	top=0pt,
	bottom=0pt,
	boxsep=3pt,
	colback=green!3!white,
	frame hidden,
	before skip=10pt plus 2pt,after skip=10pt plus 2pt,
	#1
}
\makeatother

\makeatletter
\newcommand\tabfill[1]{%
	\dimen@\linewidth
	\advance\dimen@\@totalleftmargin
	\advance\dimen@-\dimen\@curtab
	\parbox[t]\dimen@{#1\ifhmode\strut\fi}%
}
\makeatother

\usepackage{tabularx}
\usepackage{imakeidx}
\usepackage{hyperref}
\hypersetup{colorlinks={true},linkcolor={blue},citecolor={blue},urlcolor={blue}, linktoc=all}  
 
 \usepackage[capitalize, nameinlink]{cleveref}
 \crefdefaultlabelformat{#2\textbf{#1}#3}
 \crefname{figure}{Figure}{Figures} 
 \Crefname{figure}{Figure}{Figures}
 \crefname{table}{Table}{Tables}
 \Crefname{table}{Table}{Tables}
 \crefname{section}{\S\hspace{-1mm}}{\S\hspace{-1mm}}
 \Crefname{section}{\S\hspace{-1mm}}{\S\hspace{-1mm}}
 \crefname{equation}{}{}
 \Crefname{equation}{}{}
 \crefname{example}{Geometric Pattern}{Geometric Patterns} 
 \Crefname{example}{Geometric Pattern}{Geometric Patterns}

 \usepackage{graphicx,tipa}

 \usepackage{footnote}

\begin{document}

\title{\textbf{Volumes of Solid Objects in Elamite Mathematics}}

\author{Nasser Heydari\footnote{Email: nasser.heydari@mun.ca}~ and  Kazuo Muroi\footnote{Email: edubakazuo@ac.auone-net.jp}}

\maketitle

\begin{abstract}
 This article studies   three-dimensional objects and their volumes in Elamite mathematics,  particularly  those found in the Susa Mathematical Tablet No.\,14 (\textbf{SMT No.\,14}). In our discussion, we identify some basic solids whose volumes have been correctly computed in Babylonian and Elamite mathematics.  We also show that the Elamite scribes knew the right formula for calculating the volume of a certain pyramid which is a rare phenomenon occurring in the Babylonian  mathematical tablets. 
\end{abstract}

\section{Introduction}
\textbf{SMT No.\,14}   is one of   26 clay tablets excavated from Susa in  southwest Iran by French archaeologists in 1933. The texts of all the Susa mathematical tablets (\textbf{SMT}) along with their interpretations were first published in 1961 (see \cite{BR61}).

This tablet\footnote{The reader can see the new  photos of this tablet on the website of the Louvre's collection. Please see \url{https://collections.louvre.fr/en/ark:/53355/cl010186539} for obverse  and   reverse.}     consists of two problems both of which concern the volume of  an imaginary large grain-heap, whose length, width, and height are known. From a  mathematical point of view, the first problem is very important because the volume of a pyramid  is correctly calculated. Unfortunately, the second problem is nearly unintelligible to us because most parts of the text are lost.

\section{Volumes of Solid Objects}
\begin{figure}[H]
	\centering
	\includegraphics[scale=1]{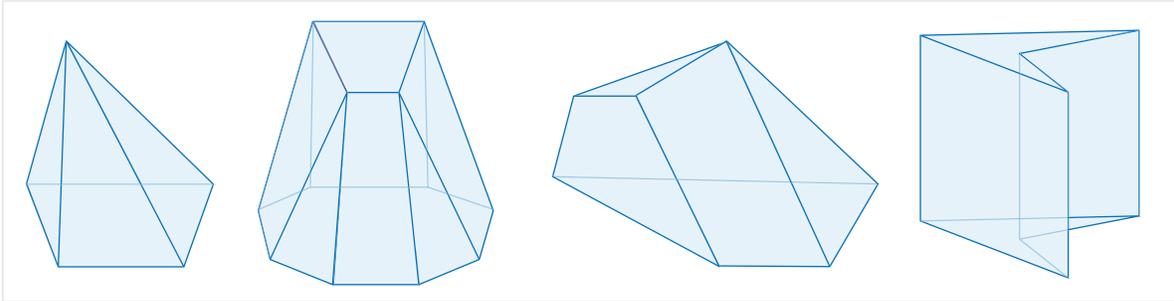}
	\caption{Examples of polyhedra}
	\label{Figure1}
\end{figure}

Usually, a \textit{solid} is any limited portion of   three-dimensional space bounded by surfaces. If all boundary surfaces forming  a solid are planes, it  is called a \textit{polyhedron}.\footnote{For a reference on polyhedra, the interested reader can consult \cite{Cro99}.} The intersections of planes in a polyhedron are called \textit{edges} and the polygons formed by the edges are \textit{faces}. The intersections of the edges are the \textit{vertices} of the polyhedron. The angle between two faces of a polyhedron is a \textit{dihedral angle}. The space near a vertex of a polyhedron is called a \textit{solid angle} or \textit{polyhedral angle}, which is formed by the intersection of at least three faces. A polyhedron is convex if the line segment connecting any pair of its points is contained completely within the polyhedron. Here, we mostly consider convex polyhedra.   \cref{Figure1} shows four polyhedra three of which are convex but the one on the right is not.
 
There is an interesting relation between the numbers of faces  $f$, edges $e$ and vertices $v$ of a polyhedron known as the  \textit{polyhedron formula}, observed by the Swiss mathematician Leonhard Euler  (1707--1783), saying that 
\begin{equation}\label{equ-SMT14-aaaa}
	v-e+f=2. 	
\end{equation}
Note that  for the  second polyhedron in \cref{Figure1}, we have $ 12-20+10= 22-20=2$ satisfying the Euler's polyhedron formula \cref{equ-SMT14-aaaa}.

\begin{figure}[H]
	\centering
	\includegraphics[scale=1]{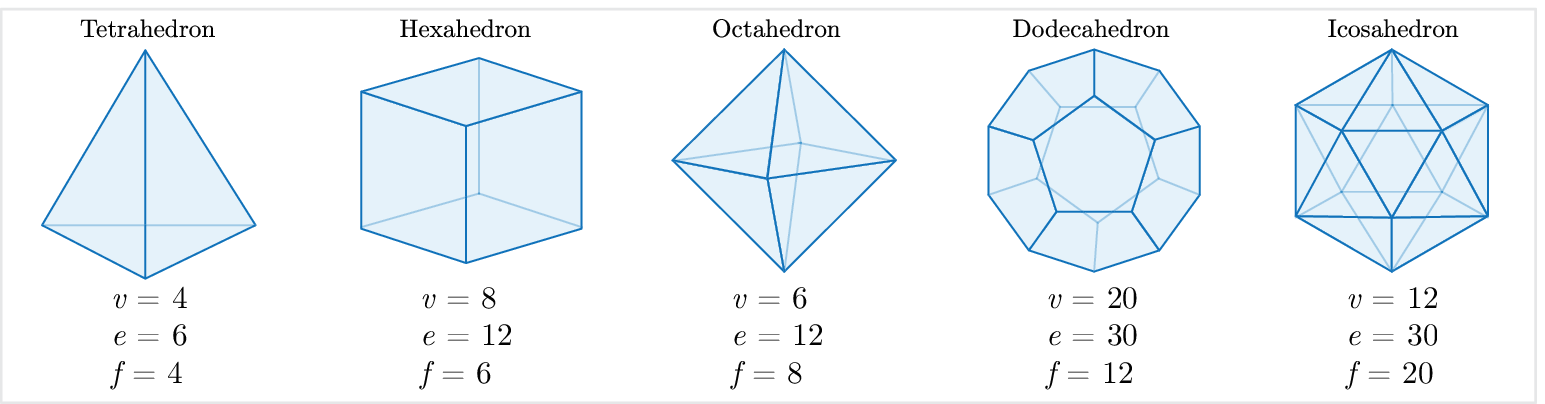}
	\caption{Platonic solids}
	\label{Figure2}
\end{figure}

Similar to regular polygons, we can define \textit{regular polyhedra} as those polyhedra all of whose faces are congruent regular polygons and all of whose polyhedral angles are equal. There is a   classic theorem in geometry stating that there are only five regular  convex polyhedra: tetrahedron with 4 faces, hexahedron (cube) with 6 faces, octahedron with 8 faces, dodecahedron  with 12 faces and  icosahedron with 20 faces (see \cref{Figure2}).

These special polyhedra are usually   called   \textit{Platonic solids}. These five regular polyhedra are named after the Greek philosopher Plato (circa 428--348 BC) who was so captivated by their perfect forms   that in his dialogue \textit{Timaeus} he associates them with what, at that time, were believed to be the basic elements: earth, fire, air, water and ether. Recently, there has been an increasing interest in a class of polyhedra  (including the Platonic solids)  because they arise naturally in a number of diverse problems in physics, chemistry, biology as well as a variety of other disciplines (see \cite{AS03}). Although these polyhedra bear the name of Plato, the Greek mathematician Theaetetus of Athens (circa 417--369 BC) was the first to give a mathematical description of all five regular polyhedra and may also have been responsible for the first known proof that no other  regular polyhedra exist (see \cite{All89,Hea21}). For a proof of this theorem, see \cite{Har00}.

The \textit{volume} of a solid is a quantity for the amount of space occupied by it which is usually measured as the number of times it contains a chosen solid as the unit of volume.   Nowadays, volumes are mostly  expressed in cubic units such as $cm^3$, $m^3$ and so on depending on the  length unit one considers for the three dimensions  of the unit cube.

A \textit{cuboid} is a polyhedron  with six faces all of which  are rectangles and all of whose dihedral angles are right angles. A  cuboid has three main dimensions, i.e.,  the length, the width and the height.     The volume of a cuboid  $\Omega$ with dimensions $a,b,c$ is simply given by
$$ V_{\Omega}=abc.$$ 
This is due to the simple observation that a cuboid with integer  dimensions $m,n,p$ contains $m\times n\times p$ numbers of the unit cube.

\begin{figure}[H]
	\centering
	\includegraphics[scale=1]{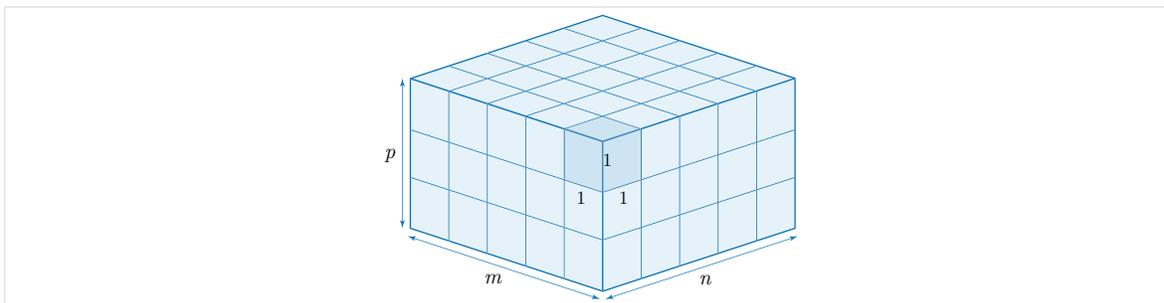}
	\caption{Volume of a cuboid}
	\label{Figure3}
\end{figure} 
 
\begin{figure}[H]
	\centering
	\includegraphics[scale=1]{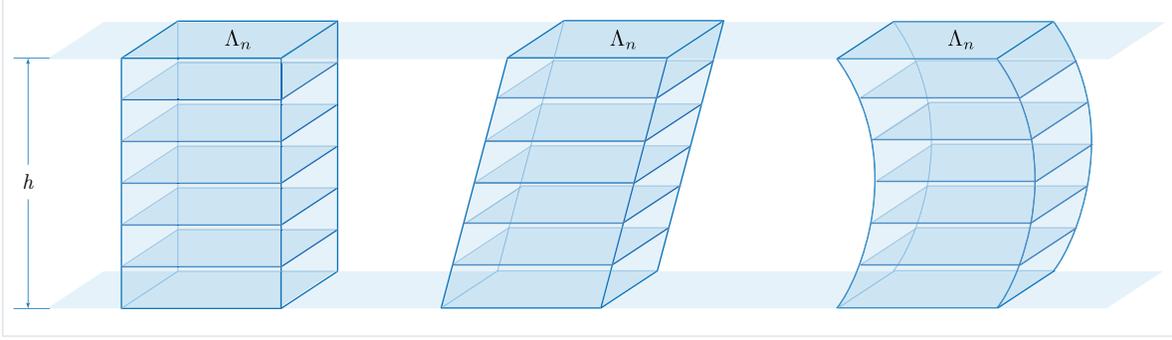}
	\caption{Cavalieri's principle}
	\label{Figure4}
\end{figure}

Although the volume of a cuboid follows directly from the definition, to find the volumes of other solids,  one needs to apply    \textit{Cavalieri' principle} attributed to the Italian mathematicians Bonaventura Cavalieri (1598--1647). This states that    if two solids are bounded between two parallel planes and if any plane parallel to the boundary planes intersects the two solids in two  cross-sections with equal areas, then the two solids have the same volume (see \cref{Figure4}). 
 
\begin{figure}[H]
 	\centering
 	\includegraphics[scale=1]{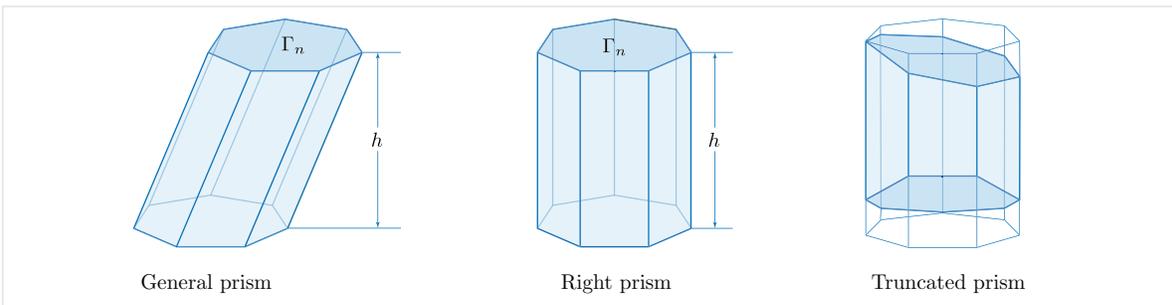}
 	\caption{Different kinds of prisms}
 	\label{Figure5}
 \end{figure}

 A \textit{prism}  $\Gamma$ is a polyhedron two of whose faces are parallel equal polygons, say $\Gamma_n$. In fact, if we    translate $ \Gamma_n$ along a fixed direction, the obtained solid is a prism. In a general prism, all faces other than the two parallel ones are parallelograms. If these parallelogram faces of a prism are rectangles, it is called  a \textit{right prism}.  The distance between the two parallel faces is called the \textit{height} of the prism denoted by $h$.  The volume of a general prism is obtained by the general rule ``area of base times height'': 
\begin{equation}\label{equ-SMT14-aaa}
	V_{\Gamma}=h\times S_{\Gamma_n}. 	
\end{equation}
This easily gives the volumes of special  prisms such as triangular prisms,   parallelepipeds and cuboids.  A \textit{truncated prism}  is  obtained from a prism by cutting it with two nonparallel planes  as shown in  \cref{Figure5}. 
 
\begin{figure}[H]
	\centering
	\includegraphics[scale=1]{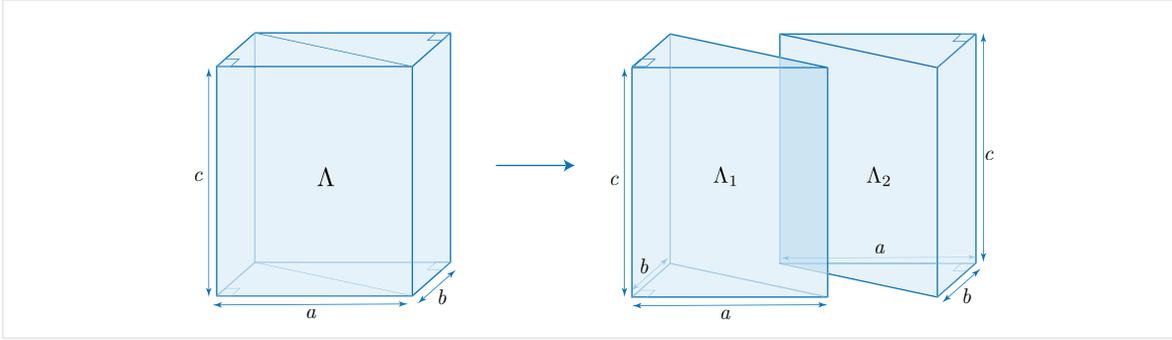}
	\caption{Cutting a cuboid into two equal right triangular prisms}
	\label{Figure6}
\end{figure}

The volume formula for a prism can be proved by using the fact that any   cuboid $\Lambda$ with sides $a,b,c$ can be diagonally  cut off into two equal triangular right prisms, say $\Lambda_1$ and $\Lambda_2$, whose bases are right triangles with sides $a$ and $b$ (see \cref{Figure6}).  So, the common volume of these two triangular prisms is 
$$V_{\Lambda_1}=V_{\Lambda_{2}}=\frac{1}{2}abc = \left(\frac{1}{2}ab\right)c,$$
which is area of base  times height. Any polygon can be partitioned into triangles, so  this, together with   Cavalieri's principle, can be used to prove   the general formula \cref{equ-SMT14-aaa}.

\begin{figure}[H]
	\centering
	\includegraphics[scale=1]{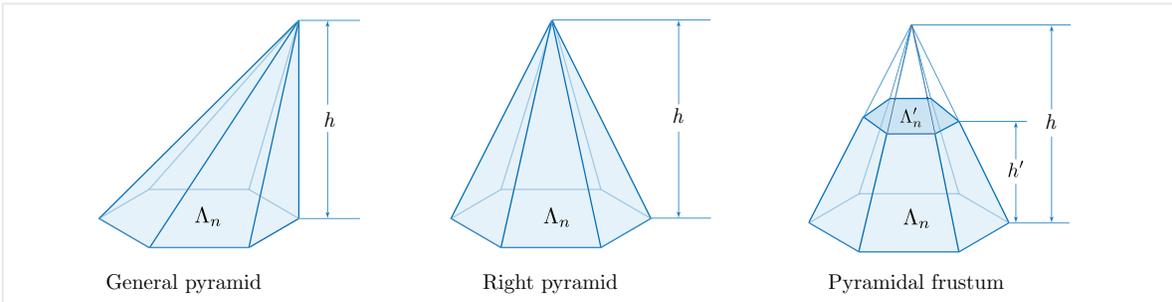}
	\caption{Different kinds of pyramids}
	\label{Figure7}
\end{figure}

A \textit{pyramid} $ \Lambda$ is a polyhedron obtained by connecting a point $p$ outside a polygon $ \Lambda_n$   by straight lines to all  of its points. This point $p$ is the \textit{apex} and the polygon $\Lambda_n$ is the base of the pyramid. The distance between the apex and the base is the height of the pyramid. Note that all faces of a pyramid containing the apex are triangles. In  a right pyramid,  the line connecting the apex and the centroid of the base is perpendicular to the base.  A \textit{pyramidal frustum} is obtained from a pyramid by cutting it with a plane parallel to its base (see \cref{Figure7}). 

\begin{figure}[H]
	\centering
	\includegraphics[scale=1]{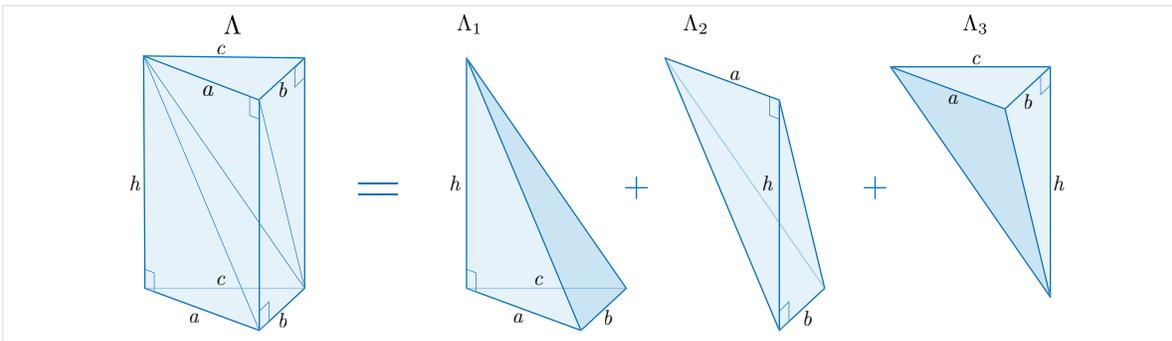}
	\caption{Dividing a prism into three  pyramids with equal volumes}
	\label{Figure8}
\end{figure}

One can make the observation that any right triangular prism $\Lambda$  is the union of three right triangular pyramids $\Lambda_1,\Lambda_2,\Lambda_3 $ with equal volumes as shown in \cref{Figure8}.

In fact, on the one hand, $\Lambda_1$ and $\Lambda_3$ have the same height $h$ and equal triangular base with sides $a,b,c$, and on the other hand, $\Lambda_2$ and $\Lambda_3$ have the same height $a$ and the same  triangular base with height $b$ and base $h$.  So by   Cavalieri's principle $V_{\Lambda_1}=V_{\Lambda_2}=V_{\Lambda_3} $. Now, since any pyramid can be divided into a finite number of triangular pyramids by triangulation of its polygonal base, one can use the previous fact to show that the volume of a general pyramid $ \Lambda$ with polygonal base $\Lambda_n$ and height $h$   is 
$$V_{\Lambda}=\frac{1}{3}\times h \times S_{\Lambda_n}.$$

\begin{figure}[H]
	\centering
	\includegraphics[scale=1]{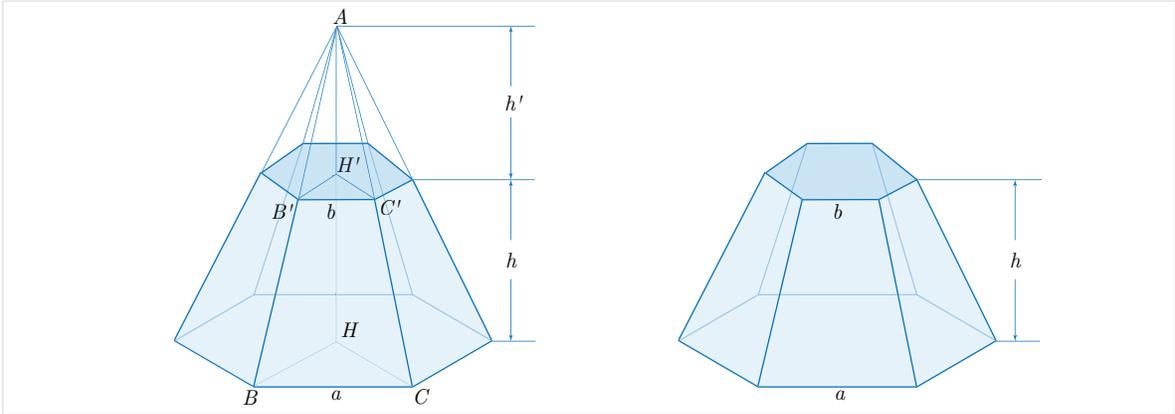}
	\caption{Volume of a pyramidal frustum whose base and top are regular polygons}
	\label{Figure9-2}
\end{figure}

For a pyramidal frustum whose base and top are regular $n$-gons with sides $a$ and $b$ respectively and whose height is $h$, one can easily find a volume formula with respect to $n$, $a$, $b$ and $h$. In fact, consider  \cref{Figure9-2}  in which the frustum with  $\Gamma_n(a)$ and $\Gamma_n(b)$ as the base and top is a part of pyramid with height $h+h'$ and base $\Gamma_n(a)$. Since the two   triangles $\triangle AH'B'$ and  $\triangle AHB$ in the figure are similar,   we have
\[ \frac{\overline{AH'}}{\overline{H'B'}} = \frac{\overline{AH}}{\overline{HB}}.\]
Also, $\overline{AH'}=h'$, $\overline{AH}=h+h'$, $\overline{B'H'}=\frac{b}{2\sin(\frac{\pi}{n})}$, and $\overline{BH}=\frac{a}{2\sin(\frac{\pi}{n})}$, so it follows from the last equality   that
\[ h'=\frac{bh}{a-b}. \]

On the other hand, since $S_{\Gamma_n(a)}=\frac{na^2}{4}\cot(\frac{\pi}{n})$ and  $S_{\Gamma_n(b)}=\frac{nb^2}{4}\cot(\frac{\pi}{n})$, and the volume of the frustum is the volume of a pyramid with height $h$ and   base  $\Gamma_n(a)$ minus that of the small pyramid with height $h'$ and base $\Gamma_n(b)$,   we can write
\[ V=\frac{h+h'}{3}\times \frac{na^2}{4}\cot\left(\frac{\pi}{n}\right) - \frac{h'}{3}\times \frac{nb^2}{4}\cot\left(\frac{\pi}{n}\right) \]
which can be simplified as
\begin{equation}\label{equ-a}
	V = \frac{nh}{12}\cot\left(\frac{\pi}{n}\right)\left(a^2+ab+b^2\right). 
\end{equation}

\begin{figure}[H]
	\centering
	\includegraphics[scale=1]{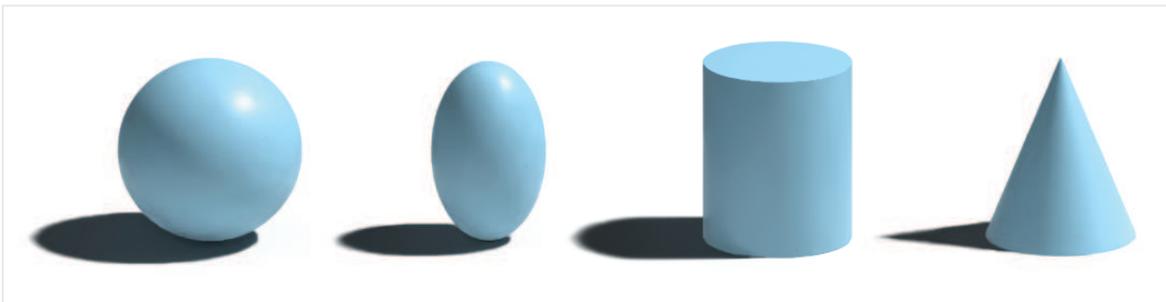}
	\caption{Different kinds of solids of rotation}
	\label{Figure9}
\end{figure}

Besides polyhedra, other solids   can be obtained by rotating   closed surfaces  around a fixed direction (the \textit{axis of rotation}) in   3-space,     generally known as \textit{solids of rotation}.   Examples of such solids of rotation are spheres, ellipsoids, cylinders and cones,  obtained by rotating a semicircle, a semi-ellipse, a rectangle and a right triangle respectively around a fixed direction as shown in \cref{Figure9}. If the radius of the   semicircle is $r$, the volume of generated sphere is $\frac{4}{3}\pi r^3$. It was Archimedes who showed that the volume of a sphere is  equal to twice the volume between the sphere and its circumscribed cylinder.

The volume of the cylinder obtained by rotating a rectangle with sides $r$ and $h$ around the latter  is $\pi r^2 h$ which follows from  Cavalieri's principle. If the base and the height of a right triangle are $r$ and $h$ and we rotate it around its height, the volume of the generated cone is  $\frac{1}{3}\pi r^2 h $ which is similar to the volume formula of a  pyramid. One way to intuitively  see why this is true is to place two containers in the shapes of a cone and a cylinder with the same height and circular base side by side. By filling the cone with water and then pouring it into the cylinder, it is apparent that the level of water   reaches   one-third of the height of the cylinder. This suggests that the volume of a cone with height $h$ and radius $r$ is  that of a cylinder with height $\frac{h}{3}$  and the same  radius $r$.

\section{Volumes in ancient mathematics}
Although the above-mentioned volume formulas for basic solids can be   obtained by    using  elementary techniques from solid geometry\footnote{The interested reader can consult \cite{KB38} for a detailed discussion on this topic.} or   double   integration from calculus, their   appearance  in mathematics can be traced back to Greek mathematicians such as  Euclid (in Book XII of the \textit{Elements}) and  Eudoxus of Cnidus (circa 395--337 BC)     who  developed solid geometry by establishing the usual proportions and volume relations for solid figures. For example, Eudoxus  proved that the volume of a sphere is proportional to the cube of its radius and the volume of a pyramid is the one-third of that of a prism with  the same base and height.

Whereas the Greek mathematicians conducted a detailed study of solids and their volume,   there is  evidence  that the Babylonians and Egyptians used certain formulas concerning  the volume  of cubes,  prisms and cones   long before the Greeks (see \cite{Fri05-1}).  Besides the purely mathematical point of view,   some of these formulas seemed to be  used in calculations regarding construction projects such as digging a hole or a canal or building a brick wall. Such   volume formulas would have been of great advantage to builders. For example, by using the dimensions of a cubical  brick and an intended wall, they could estimate the total number of bricks needed to finish the wall. In addition, such        data could be used to determine the number of workers required for any project.

\begin{figure}[H]
	\centering
	\includegraphics[scale=1]{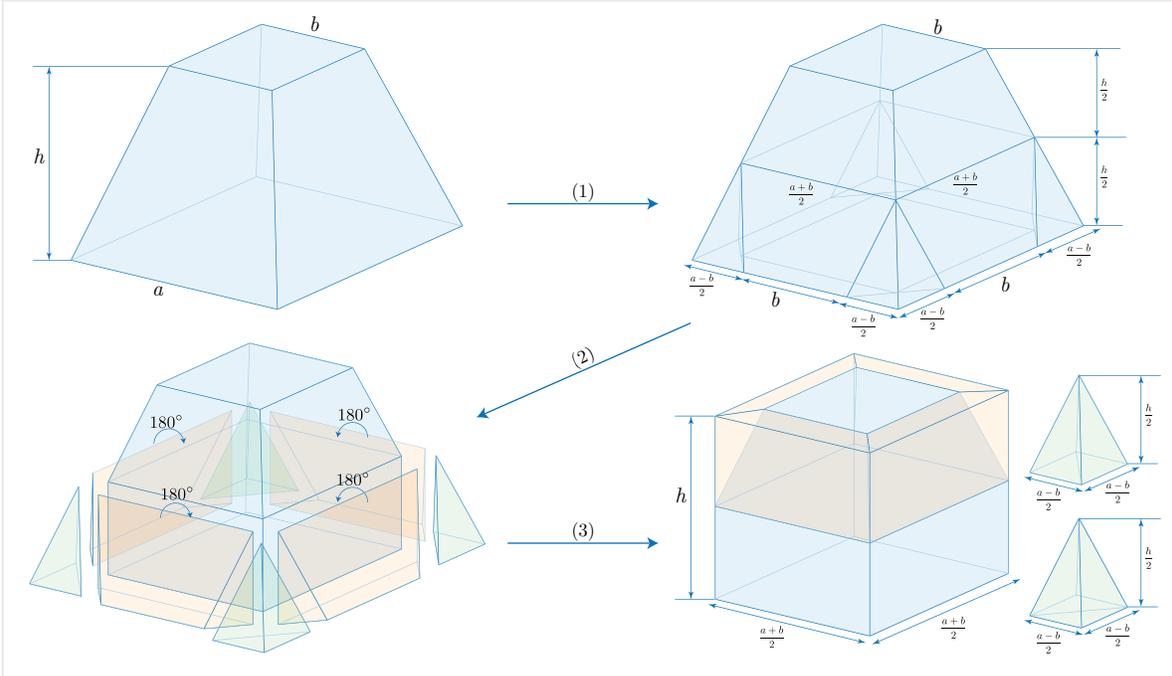}
	\caption{Babylonian formula for the volume of a pyramidal frustum}
	\label{Figure9-3}
\end{figure}

Unsurprisingly, calculating the volume of pyramids was a challenge for ancient scribes. Apparently, the  Egyptians and the Babylonians were aware of  the   volume formula  for a pyramid but they did not use it explicitly. One can see this implicit usage in the Babylonian and the Egyptian formulas for a pyramidal frustum, say $\Delta$, with square base of side $a$ and square top of side $b$ and height $h$ (see \cref{Figure9-3}). As is known from mathematical tablet \textbf{BM 85194}\footnote{For an interpretation of a part   of this text, see the appendix section.}, the  Babylonians   used the  volume formula 
\begin{equation}\label{equ-ab}
	V^B_{\Delta}=\left[\left(\frac{a+b}{2}\right)^2+\frac{1}{3}\left(\frac{a-b}{2}\right)^2\right]h
\end{equation}
 for the volume of the pyramidal frustum $\Delta$. 

One way to obtain this formula is the ``cut and paste'' method given in \cref{Figure9-3}.\footnote{For more details, see \cite[Vol. II, pp. 332-334]{Hea21};\cite[pp. 88-93]{Mur00};  or \cite[pp. 136-141]{Vai61}.} A pyramidal frustum is first cut vertically at its four corners to get four pyramids  with height $\frac{h}{2}$ and square bases  with side $\frac{a-b}{2}$. Next, we slice off   four extra parts on four faces of the lower half frustum to get four truncated prisms. Then these  four polyhedra are rotated $180^{\circ}$ around their top edges and   attached to the faces of the upper half frustum as shown in the figure.  We can also attach   two  small triangular pyramids on their proper faces to get a pair of pyramids with square bases. At the end, we obtain a cube with dimensions  $h$, $\frac{a+b}{2}$ and   $\frac{a+b}{2}$ plus two pyramids of square bases with sides $\frac{a-b}{2}$ and height $\frac{h}{2}$. The total volume of these solids is exactly the value given in \cref{equ-ab}. The second part of the formula clearly proves that the Babylonians knew the formula of a pyramid with a square base.

On the other hand,  according to the \textit{Moscow papyrus}\footnote{The Moscow Mathematical Papyrus  is an ancient Egyptian mathematical papyrus containing several problems in arithmetic, geometry, and algebra. It is held in the collection of the Pushkin State Museum of Fine Arts in Moscow.}, the Egyptians applied the formula
\begin{equation}\label{equ-ac}
V^E_{\Delta}=\frac{h}{3}\left(a^2+ab+b^2\right)
\end{equation}
for the same pyramidal frustum (note that this is obtained from \cref{equ-a} by setting $n=4$). Although  the formula is different than that of the Babylonians, the two formulas are   equivalent.  We can use the ``cut and paste'' method shown in   \cref{Figure9-3} to obtain the formula (see \cite[pp. 35-39]{NM14}). The only thing in this method is that we replace four pyramids with heights $h$ in the corners with four cubes with the same base but height $\frac{h}{3}$. The other steps are just cutting and pasting processes as shown in the figure.

\begin{figure}[H]
	\centering
	\includegraphics[scale=1]{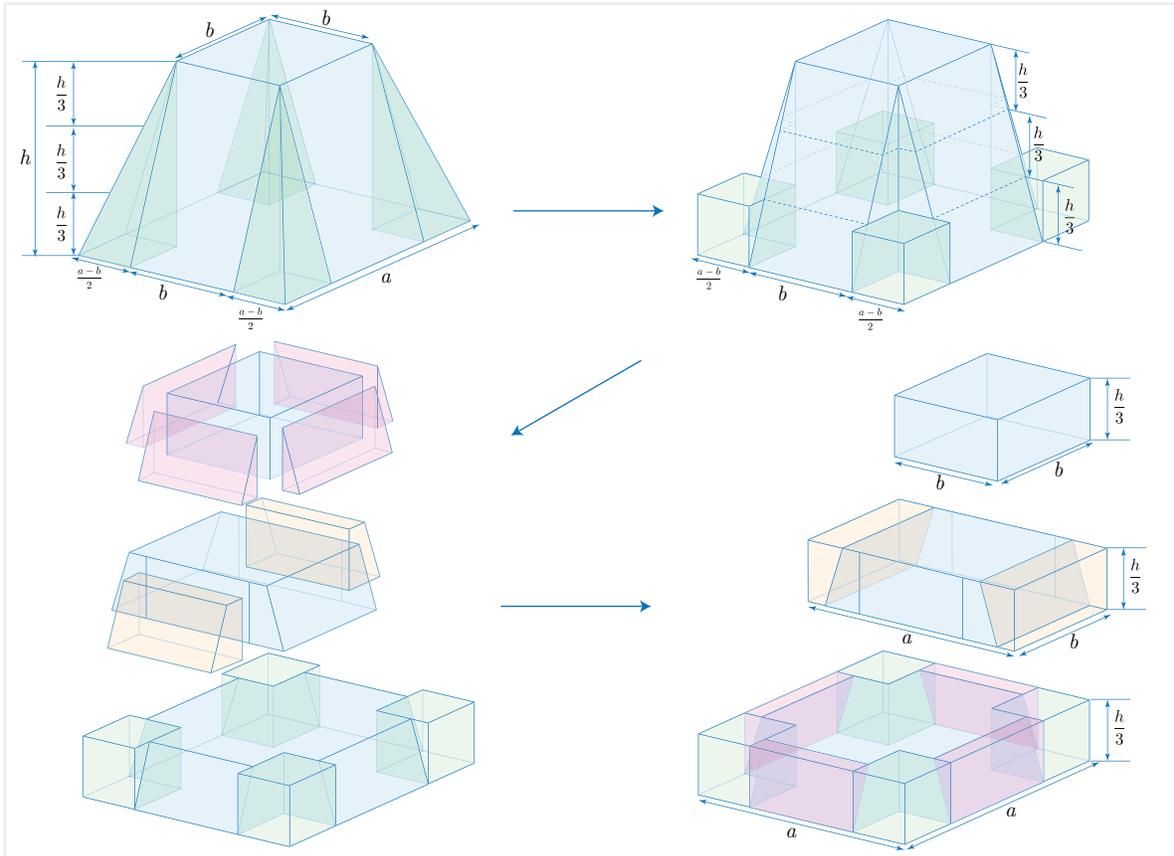}
	\caption{Egyptian formula for the volume of a pyramidal frustum}
	\label{Figure9-1}
\end{figure}

\begin{figure}[H]
	\centering
	\includegraphics[scale=1]{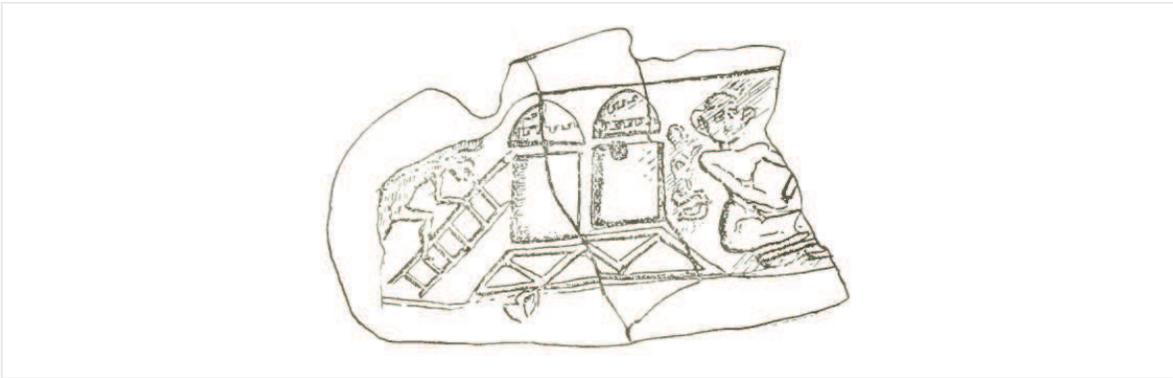}
	\caption{A pair of granaries on a Susa seal  dated to circa 3500-3000 BC}
	\label{Figure10-1}
\end{figure}

  Although  no explicit formula or direct calculation for the volume of a pyramid is  given, there are clay tablets that address problems regarding the volumes of solids whose calculations involve the volume of a pyramid.  Two   such tablets are   \textbf{BM 96954} and \textbf{SMT No.\,14}    dealing  with   the volume of   granaries.  In both mathematical and nonmathematical texts, different shapes are suggested for ancient   granaries. The common   shape of a cylinder with a domed top is one of the oldest and  can be found on  seal imprints. For example, in figure 222 of \textbf{MDP XVI} (see \cite[Plate 15]{Leg21}),   a worker is climbing up a ladder to put  grain into a pair of cylindrical  granaries (see \cref{Figure10-1}). This ancient seal  from Susa caused some   scholars to suggest that the cylinder was the likely shape of granaries in ancient Elam. Unfortunately, this   led them to misunderstand some mathematical texts and so to   mistakes in their interpretations (see \cite[TEXTE XIV]{BR61}, for example).

\begin{figure}[H]
	\centering
	\includegraphics[scale=1]{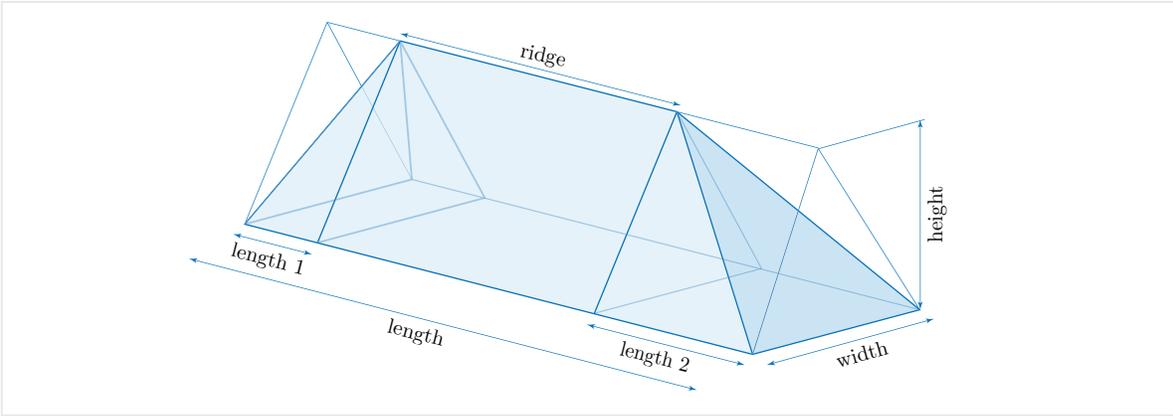}
	\caption{A truncated  triangular prism and its dimensions}
	\label{Figure10}
\end{figure}

In this process of computing the volume of a truncated triangular prism (see  \cref{Figure10}),   the Babylonian and Elamite scribes usually  needed  to calculate the volume  of a special  rectangular pyramid.  This   truncated triangular prism  was a   solid shape familiar  to the Elamites and the Babylonians and might have served   as a pattern for building storage facilities such as granaries  due to the rigid  structure provided by  this shape.

Let us compute the volume of   the truncated triangular prism in \cref{Figure10} by  assuming that its height is $h$, its width is $y$, its ridge is $x$, and  its length is $z$.  We also denote the other two parts of the length by   $x_1$ and $x_2$ (that is, the lengths of two right rectangular pyramids attached to the middle triangular prism). Clearly, $z=x+x_1+x_2$.  In general, $ x_1$ and $x_2$ may not be equal although Babylonian and Elamite scribes considered the special case  $x_1=x_2$.   First, note that our solid consists of three parts: the middle part is a triangular prism and the left and right parts are two   rectangular pyramids. The volume of the triangular prism is $V_0=\frac{xyh}{2}$ while those of the two rectangular pyramids are $V_1=\frac{1}{3}\times  x_1\times y \times h=\frac{x_1yh}{3} $ and similarly $V_2=\frac{x_2yh}{3}  $. Therefore,  
\[ V=V_0+V_1+V_2=\frac{xyh}{2}+\frac{x_1yh}{3}+\frac{x_2yh}{3}=\frac{yh}{6}\times (3x+2x_1+2x_2) \]
implying that
\begin{equation}\label{equ-SMT14-aa}
	V=\frac{2zyh+xyh}{6}.	
\end{equation}

As we see in the next section, the Susa scribes   correctly computed the volume of such a truncated triangular prism in which the two pyramids are equal (that is, $x_1=x_2$). Further, mathematical analysis of the text   shows that the Susa scribes    knew  the  formula for the volume of a specific  rectangular pyramid whose height and   sides of the   base are all the same number $h$. They   used the formula $\frac{h^3}{3} $ for the volume which is obtained from the usual formula ``one-third of height times the area of base''.  

\section{\textbf{SMT No.\,14}}
As mentioned, this text    consists of two problems both of which concern the volume\index{volume of a grain-heap} of an imaginary   grain-heap\index{grain-heap}, whose length\index{length}, width\index{width}, and height\index{height} are 10 {\fontfamily{qpl}\selectfont nindan}\index{nindan (length unit)} ($\approx 60m $), 6 {\fontfamily{qpl}\selectfont nindan}\index{nindan (length unit)} ($ \approx 36m$), and 3 {\fontfamily{qpl}\selectfont nindan}\index{nindan (length unit)} ($ \approx 18m$) respectively. From the mathematical point of view, the first problem is very important because the volume\index{volume of a pyramid} of a pyramid\index{rectangular pyramid} is correctly calculated and the mathematical technical term {\fontfamily{qpl}\selectfont a-r\'{a} \textit{kayyam\={a}num}}, thought to be the Akkadian translation of the Sumerian word\index{Sumerian word} {\fontfamily{qpl}\selectfont a-r\'{a}-gub-ba}, occurs in a complete form. The second problem is sadly nearly unintelligible to us as most parts of the text are lost.

\subsection{Transliteration}\label{SS-TI-SMT14}

\begin{note1} 
	\underline{Obverse:  Lines 1-18}\\
	(L1)\hspace{2mm} gur$_7  $ \textit{a-na} 14,24 sahar 3 gi(sic) \textit{me-l}[\textit{a-a-am} gar]\\
	(L2)\hspace{2mm} \textit{a-na} 14,24 sahar u\v{s} sag \textit{\`{u} ka-aq-qa-da}\\
	(L3)\hspace{2mm} \textit{mi-na} gar za-e igi-12 \textit{\v{s}u-up-li pu-\c{t}\'{u}-\'{u}r}\\
	(L4)\hspace{2mm} 5 \textit{ta-mar} 5 \textit{a-na} 14,24 sahar \textit{i-\v{s}\'{i}-ma}\\
	(L5)\hspace{2mm} 1,12 \textit{ta-mar} 3 gi(sic) \textit{me-la-a-am} nigin 9 \textit{ta-mar}\\
	(L6)\hspace{2mm} 9 \textit{a-na} 3 \textit{me-le-e te-er-ma} 27 \textit{ta-mar}\\
	(L7)\hspace{2mm} \textit{i-na} 1 a-r\'{a} \textit{ka-aiia-ma-ni} 20 [sa]har \textit{\v{s}\`{a}-lu-u\v{s}-ti}\\
	(L8)\hspace{2mm}\label{SMT14-obv-L8} \textit{\v{s}\`{a} ka-aiia-ma-ni} \'{a} \textit{tu-u\c{s}-\c{s}a-bu} zi\\
	(L9)\hspace{2mm} 40 \textit{ta-mar} 40 \textit{a\v{s}-\v{s}um} 2 sig gur$_7  $ \textit{a-na} 2 tab-ba\\
	(L10)\hspace{0mm} 2(sic),20 \textit{ta-mar} 1,20 \textit{a-na} 27 \textit{i-\v{s}\'{i}-ma} 36 \textit{ta-mar}\\
	(L11)\hspace{0mm} 36 \textit{i-n}[\textit{a}] 1,12 zi 36 [\textit{ta-mar tu}]-\textit{\'{u}r-ma}\\
	(L12)\hspace{0mm} 3 \textit{me-la-a-am} nigin 9 \textit{ta}-[\textit{mar} igi-9 \textit{pu-\c{t}}]\textit{\'{u}-\'{u}r}\\
	(L13)\hspace{0mm} 6,40 \textit{ta-mar} 6,40 \textit{a-na} 36 \textit{i-\v{s}}[\textit{\'{i}-ma}]\\
	(L14)\hspace{0mm} 4 \textit{ta-mar} 4 \textit{ka-aq-qa-du} 3 \textit{me-l}[\textit{a-a-am} gar]\\
	(L15)\hspace{0mm} [\textit{a\v{s}-\v{s}u}]\textit{m} \textit{i-na am-ma-at am-ma-a}[\textit{t i-ku-lu}]\\
	(L16)\hspace{0mm} [\textit{a-na} 2] tab-ba 6 \textit{ta-mar} 6 sag 6 [\textit{a-na} 4]\\
	(L17)\hspace{0mm} [\textit{ka-aq-qa}]-\textit{di} dah 10 \textit{ta-mar} 10 [u\v{s}]\\
	(L18)\hspace{0mm} [$ \cdots $ $ \cdots $ $ \cdots $ $ \cdots $]\\
	\underline{Reverse:  Lines 1-7}\\
	(L1)\hspace{2mm} [$ \cdots $ $ \cdots $ T]A(?) [$ \cdots $ $ \cdots $ $ \cdots $]\\
	(L2)\hspace{2mm} [$ \cdots $ $ \cdots $] 3 \textit{me-la-a-am} [\textit{a-na} 12]\\
	(L3)\hspace{2mm} [\textit{i-\v{s}\'{i}-ma}] 36 \textit{ta-mar} 3[6 \textit{a-na} 24 \textit{i-\v{s}\'{i}-ma}]\\
	(L4)\hspace{2mm} [14],24 \textit{ta-mar} sahar 14,2[4 $ \cdots $ $ \cdots $]\\
	(L5)\hspace{2mm}  \textit{a-na} 8 \textit{na-a\v{s}-pa-ak} gur$_7  $ [\textit{i-\v{s}\'{i}-ma}]\\
	(L6)\hspace{2mm} 1,55,12 gub-\textit{ma} 20,30(sic) [gur$_7  $]\\
	(L7)\hspace{2mm} \textit{\`{u}} 2-\textit{\v{s}u} 24 gur \textit{\v{s}e-u}[\textit{m ki-a}]-\textit{am ne}-[\textit{p\'{e}-\v{s}um}]
	
\end{note1}

\subsection{Translation}\label{SS-TR-SMT14}

\underline{Obverse:  Lines 1-18}
\begin{tabbing}
	\hspace{18mm} \= \kill 
	(L1)\> \tabfill{The grain-heap. For the volume 14,24, I put down 3 ({\fontfamily{qpl}\selectfont nindan}\index{nindan (length unit)}, that is, 6) {\fontfamily{qpl}\selectfont gi} as the height.}\index{height}\index{volume of a pyramid}\index{grain-heap}\index{gi (length unit)}\\
	(L2-3)\> \tabfill{For the volume 14,24, what did I put down as the length, width, and the top? You, make the reciprocal of the (constant) 12 of depth, (and)}\index{width}\index{length}\index{reciprocal of a number}\index{depth (of a canal)}\index{volume of a pyramid} \\ 
	(L4)\> \tabfill{you see 0;5. Multiply 0;5 by the volume 14,24, and}\index{volume of a pyramid}\\ 
	(L5)\> \tabfill{you see 1,12. Square the height, 3 ({\fontfamily{qpl}\selectfont nindan}\index{nindan (length unit)}, that is, 6) {\fontfamily{qpl}\selectfont gi}, (and) you see 9.}\index{gi (length unit)}\index{height}\\
	(L6)\> \tabfill{Multiply 9 by the height 3 again, and you see 27.}\index{height}\\ 
	(L7-8)\> \tabfill{From the regular number 1, subtract 0;20 of the volume, one third of the regular (number 1) of the wing that you add (to the middle part), (and)}\index{volume of a pyramid}\\
	(L9)\> \tabfill{you see 0;40. Since there are two dilapidated parts of the grain-heap, double 0;40, (and)}\index{grain-heap}\index{wing of a grain-heap}\index{dilapidated part of a grain-heap}\\
	(L10)\> \tabfill{you see 2;20 (error for 1;20). Multiply 1;20 by 27, and you see 36.}\\ 
	(L11)\> \tabfill{Subtract 36 from 1,12, (and) you see 36. Return and}\\
	(L12)\> \tabfill{square the height 3, (and) you see 9. Make the reciprocal of 9, (and)}\index{reciprocal of a number}\index{height} \\
	(L13)\> \tabfill{you see 0;6,40. Multiply 36 by 0;6,40, (and)}\\
	(L14)\> \tabfill{you see 4. 4 is the top (length). Put down the height 3.}\index{length}\index{height}\\
	(L15-16)\> \tabfill{Since the inclination of the sides of the grain-heap is 1 {\fontfamily{qpl}\selectfont k\`{u}\v{s}} ($ \approx 50cm$) per 1 {\fontfamily{qpl}\selectfont k\`{u}\v{s}}, double the height 3, (and) you see 6. 6 is the width.}\index{height}\index{width}\index{kuz@k\`{u}\v{s} (length unit)}\index{grain-heap}\index{inclination of a plane}\\ 
	(L16-17)\> \tabfill{Add 6 to the top (length) 4, (and) you see 10. 10 is the length.}\index{length}\\
	(L18)\> \tabfill{$ \cdots $ $ \cdots $ $ \cdots $ $ \cdots $.} 
\end{tabbing}\index{normal number}\index{regular number} 

\noindent  
\underline{Reverse:  Lines 1-7}
\begin{tabbing}
	\hspace{14mm} \= \kill 
	(L1)\> \tabfill{$ \cdots $ $ \cdots $ $ \cdots $ $ \cdots $.}\\
	(L2-3)\> \tabfill{Multiply the height 3 by 12 $ \cdots $, and you see 36.}\index{height}\\
	(L3-4)\> \tabfill{Multiply 36 by 24, and you see 14,24.}\\
	(L4-7)\> \tabfill{Multiply the volume 14,24, $ \cdots $, by 8,0,0, the storage (constant) of the grain-heap, and 1,55,12,0,0 ({\fontfamily{qpl}\selectfont s\`{i}la}\index{szila@s\`{i}la (capacity unit)}) is confirmed, and 20,30 (error for 23) {\fontfamily{qpl}\selectfont gur$_7$} and 2,24 {\fontfamily{qpl}\selectfont gur}\index{gur (capacity unit)} is the barley\index{barley}. Such is the procedure.}\index{volume of a pyramid}\index{grain-heap}  
\end{tabbing}

\subsection{Technical Terms in \textbf{SMT No.\,14}} 
Before discussing the mathematical meaning of the problems, we    explain a few technical terms that occur in the text.

\subsubsection*{The length\index{length} unit {\fontfamily{qpl}\selectfont gi}\index{gi (length unit)}}
Since 1 {\fontfamily{qpl}\selectfont gi}\index{gi (length unit)}  is equal to half a {\fontfamily{qpl}\selectfont nindan}\index{nindan (length unit)} ($ \approx $ 3m), the value 3 {\fontfamily{qpl}\selectfont gi}\index{gi (length unit)} in lines 1 and 5  seems to be a mistake if it means ``3 ({\fontfamily{qpl}\selectfont nindan}\index{nindan (length unit)}, that is, 6) {\fontfamily{qpl}\selectfont gi}\index{gi (length unit)}''. However, we know that the same expressions were occasionally used in mathematical texts, for example:

\noindent
(i) 30 {\fontfamily{qpl}\selectfont gi}\index{gi (length unit)} ``0;30 ({\fontfamily{qpl}\selectfont nindan}\index{nindan (length unit)}, that is, 1) {\fontfamily{qpl}\selectfont gi}\index{gi (length unit)}'',\\ 
(ii) 30 {\fontfamily{qpl}\selectfont k\`{u}\v{s}}\index{kuz@k\`{u}\v{s} (length unit)} ``0;30 ({\fontfamily{qpl}\selectfont nindan}\index{nindan (length unit)}, that is, 6) {\fontfamily{qpl}\selectfont k\`{u}\v{s}}\index{kuz@k\`{u}\v{s} (length unit)}''. 

\noindent
Note that the basic length unit {\fontfamily{qpl}\selectfont nindan}\index{nindan (length unit)} is always omitted. Due to the absence of the ``sexagesimal\index{sexagesimal number system} point''\index{sexagesimal point} in Babylonian mathematics\index{Babylonian mathematics}, the scribe   used the clumsy expression  3 {\fontfamily{qpl}\selectfont gi}\index{gi (length unit)} in our problem. The   practice of writing a number in this way  may also be found in a Sumerian inscription written around 2400 BC (see \cite{Mur15,Mur17}).

\subsubsection*{Rectangular pyramids\index{rectangular pyramid}}
\begin{figure}[H]
	\centering
	\includegraphics[scale=1]{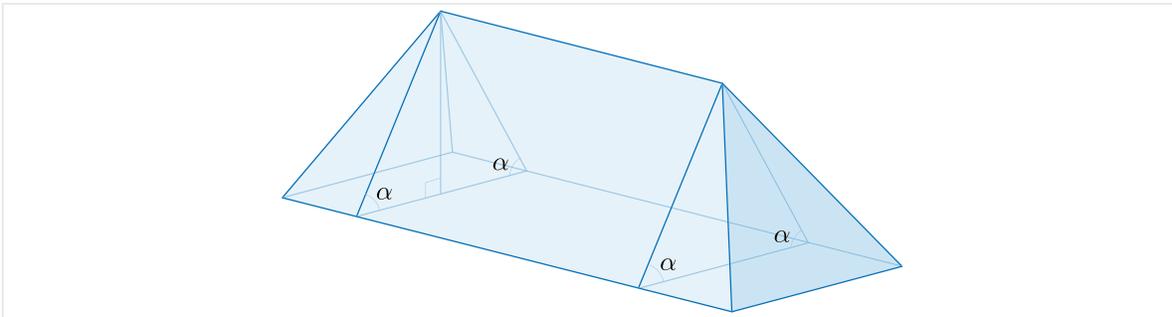}
	\caption{An imaginary grain-heap}
	\label{Figure11-1}
\end{figure}

The grain-heap\index{grain-heap} ({\fontfamily{qpl}\selectfont gur$_7 $}) consists of a right triangular prism\index{triangular prism} and  two  rectangular pyramids\index{rectangular pyramid} attached to the top and base of the triangular prism (see \cref{Figure11-1}).

In line 8, the pyramid\index{rectangular pyramid} is called {\fontfamily{qpl}\selectfont \'{a}(\textit{ahum})} ``a wing\index{wing of a grain-heap} (of the grain-heap\index{grain-heap})'' and also {\fontfamily{qpl}\selectfont sig gur$_7 $} ``a dilapidated\index{dilapidated part of a grain-heap} part of the grain-heap\index{grain-heap}'' in line 9.

\begin{figure}[H]
	\centering
	\includegraphics[scale=1]{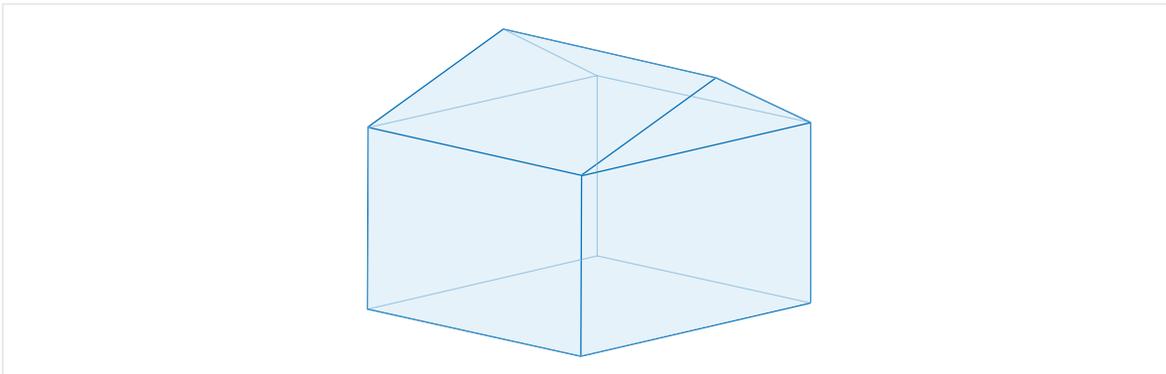}
	\caption{The  grain-heap suggested by Bruins}
	\label{Figure11-2}
\end{figure}

As mentioned, the shape in this text was misunderstood by Bruins and he mistakenly considered a solid shape as shown in \cref{Figure11-2} represented this grain-heap. Friberg believes that the imprint of a Susa seal (see \cref{Figure10-1})  was the  cause of this error   (\cite{Fri07-2}). 

\subsubsection*{Adjective {\fontfamily{qpl}\selectfont  \textit{kayyam\={a}num}}}
In several Susa mathematical texts\index{Susa mathematical texts} the Akkadian adjective {\fontfamily{qpl}\selectfont \textit{kayyam\={a}num}} ``normal\index{normal number}\index{regular number}, regular\index{normal number}\index{regular number}, usual''   occurs with the number 1 or 2 or 3. It determines the integer part of a number, namely, it is used    to refer to, for example, the number 2  in   sexagesimal\index{sexagesimal number system} numbers like  $2\times60^{\pm n}$:
\[2\times 60^{2}=2,0,0\ \ \text{and}\ \ 2\times 60^{-4}=0;0,0,0,2.\]   
Additionally, there is a possibility that 5 {\fontfamily{qpl}\selectfont a-r\'{a}} in \textbf{SMT No.\,7}\index{SMT No.g@\textbf{SMT No.\,7}} is an abbreviation for 5 {\fontfamily{qpl}\selectfont a-r\'{a} \textit{kayyam\={a}num}}, because the numbers 5 and 7 occur together as the prime\index{prime number} factors of 35, and the former is called ``the factor\index{regular factor} 5'' and the latter is not modified at all  (see  \cite{Mur92-2}). Therefore, it is highly probable that the term {\fontfamily{qpl}\selectfont a-r\'{a} \textit{kayyam\={a}num}}, which is a translation of the Sumerian term {\fontfamily{qpl}\selectfont a-r\'{a}-gub-ba}, is  specifically used for the numbers 1, 2, 3, and 5. In other words, these numbers are ``normal\index{normal number}\index{regular number}, regular\index{normal number}\index{regular number}'' in the sense that their reciprocals\index{reciprocal of a number} can be expressed by finite sexagesimal\index{sexagesimal number system} fractions\index{finite sexagesimal fraction}. The
Sumerians must have known the fact that the numbers 2, 3, and 5 are ``regular\index{normal number}\index{regular number}'' with respect to the base 60, and this  was also known  to the Susa and the Babylonian scribes.

\subsubsection*{Inclination\index{inclination of a plane} of a plane}
In line 15 occurs a technical expression of Babylonian mathematics\index{Babylonian mathematics} that defines the inclination of a plane\index{inclination of a plane}. A typical example is as follows:

\vspace{1mm}
\noindent
{\fontfamily{qpl}\selectfont \textit{i-na} 1 k\`{u}\v{s} $x$ k\`{u}\v{s} (k\'{u}) \`{i}-k\'{u}} \\
``In 1 {\fontfamily{qpl}\selectfont k\`{u}\v{s}}\index{kuz@k\`{u}\v{s} (length unit)} (in height\index{height}) it ate $x$ {\fontfamily{qpl}\selectfont k\`{u}\v{s}}\index{kuz@k\`{u}\v{s} (length unit)} (of fodder\index{fodder})''. 

\begin{figure}[H]
	\centering
	\includegraphics[scale=1]{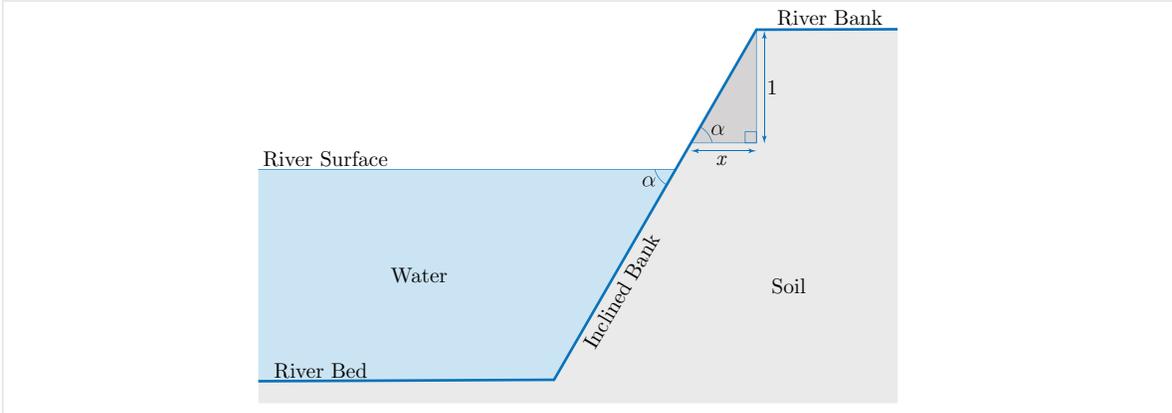}
	\caption{Erosion of a river bank}
	\label{Figure11}
\end{figure}

The underlying idea of this expression, which is also thought  to be of Sumerian
origin, is ``the water ate away the bank of a river\index{river bank}''.   It refers to the angle of erosion\index{angle of erosion}   formed between the river surface   and the incline of the   eroded bank  which is denoted by $\alpha$ in  \cref{Figure11} depicting a cross section  of both  a river and bank. 

By considering the  right triangle in  \cref{Figure11} whose bases are 1 {\fontfamily{qpl}\selectfont k\`{u}\v{s}}\index{kuz@k\`{u}\v{s} (length unit)} and $x$ {\fontfamily{qpl}\selectfont k\`{u}\v{s}}\index{kuz@k\`{u}\v{s} (length unit)}, we can compute    $ \tan(\alpha)=\frac{1}{x} $. 
In our problem since $x = 1$,   the angle between the ground and the side of the grain-heap\index{grain-heap} is  computed as $\alpha=\arctan({\frac{1}{1}})=\arctan(1)=45^{\circ}$ (see \cref{Figure11-1} and \cref{Figure12}).

\subsection{Mathematical Calculations} 
\begin{figure}[H]
	\centering
	\includegraphics[scale=1]{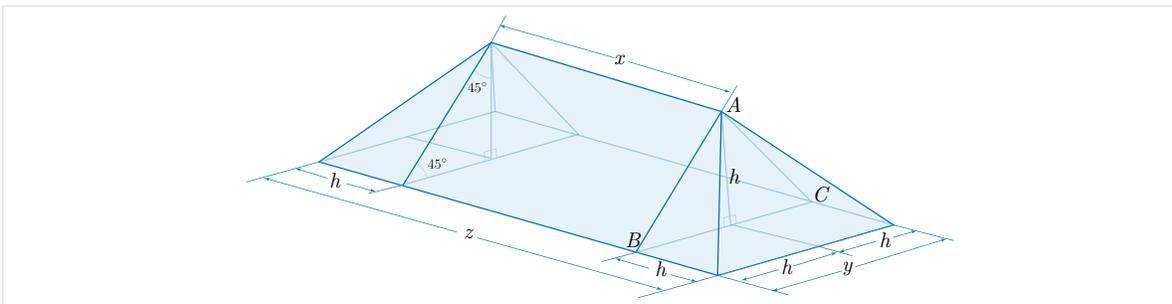}
	\caption{Dimensions of an imaginary Elamite grain-heap}
	\label{Figure12}
\end{figure}

We now analyze  the first problem of \textbf{SMT No.\,14}. In the statement of this problem only two pieces  of data are provided: the volume\index{volume of a pyramid}  and the height\index{height} of a grain-heap\index{grain-heap}, which are obviously insufficient to determine the shape of the solid figure being considered. However, careful consideration of the calculation assists in understanding the geometrical intention of the scribe. The solid figure is a truncated triangular prism  as shown in \cref{Figure12}. Note that, as we said in the previous section, both   acute angles of the right triangles in the figure are $45^{\circ}$, so they are isosceles   triangles too. 

As   in  \cref{Figure12}, we  denote the top (length\index{length}) by $x$,     the width\index{width} (of the base) by $y$,   the length\index{length} (of the base) by $z$, and  the height\index{height}  by $h$. Note that the value of height here is $h= 3$ {\fontfamily{qpl}\selectfont nindan}\index{nindan (length unit)}. Since the height $h$ bisects the base $BC$ in the isosceles right triangle $ \triangle ABC$, the scribe has rightly assumed that
\begin{equation}\label{equ-SMT14-a}
	\begin{dcases}
		y = 2h\\
		z = x + 2h.
	\end{dcases}
\end{equation}
It follows from the translation of the text that the scribe has used the following formula for the volume  $V$ of the three-dimensional solid in \cref{Figure12}:
\begin{equation}\label{equ-SMT14-ab}
	V= x h^2  + 2(1-0;20)h^3.
\end{equation}
Note that this formula is consistent with \cref{equ-SMT14-aa}, because by setting $y=2h$ and $z=x+2h$ in \cref{equ-SMT14-aa}, we get
\begin{align*}
	V&=\frac{yh}{6}(2z+x)\\
	&=\frac{(2h\times h)}{6}(3x+4h) \\
	&=xh^2+\frac{4}{3}h^3\\
	&=xh^2+2\left(1-\frac{1}{3}\right)h^3\\
	&=x h^2  + 2(1-0;20)h^3.	
\end{align*}

\begin{figure}[H]
	\centering
	\includegraphics[scale=1]{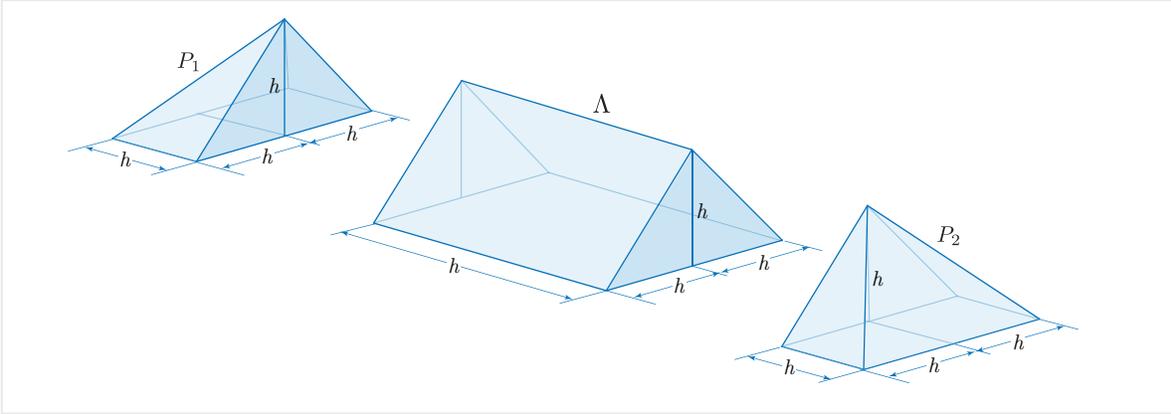}
	\caption{Splitting an imaginary grain-heap}
	\label{Figure13}
\end{figure}

Let us break down the formula \cref{equ-SMT14-ab}. First, note that our solid consists of two equal rectangular pyramids\index{rectangular pyramid}, say $P_1$ and $P_2$, as well as a triangular prism\index{triangular prism}, say $\Lambda$ (see \cref{Figure13}). 

Clearly, the volume\index{volume of a triangular prism} of $\Lambda$ is obtained by ``area of base times height'':
\begin{align*}
	V_{\Lambda}& =S_{\triangle ABC}\times x = \left[ \frac{1}{2}(h)(2h)\right]\times x 
\end{align*}    
so
\begin{equation}\label{equ-SMT14-b}
	V_{\Lambda}=h^2 x.
\end{equation}
The remaining part $ 2\left(1-\frac{1}{3}\right)h^3$ is the total volume of two equal rectangular pyramids $P_1$ and $P_2$.  So, it follows from  \cref{equ-SMT14-ab}  that 
$$V_{P_1}=V_{P_2}= (1-0;20)h^3$$
or equivalently
\begin{equation}\label{equ-SMT14-ba}
	V_{P_1}=V_{P_2}=  \left(1-\frac{1}{3}\right)h^3
\end{equation} 
which confirms that Susa scribes have computed   the correct values of these volumes.

 \subsubsection*{Formula for the Volume of a Rectangular  Pyramid}
The statement of   formula \cref{equ-SMT14-ba} suggests two facts about the volumes of   pyramids. Firstly,  the   Susa scribes have assumed that the volume of a rectangular pyramid  with height $h$, width $h$ and length $h$ is $\frac{1}{3}h^3$. Secondly, if they subtract  this value from the volume of a cube of the same dimensions, they get  the volume of a    rectangular pyramid with height $h$, length $h$ and width  $2h$. In the following, we give a geometric explanation for the first fact.    

\begin{figure}[H]
	\centering
	\includegraphics[scale=1]{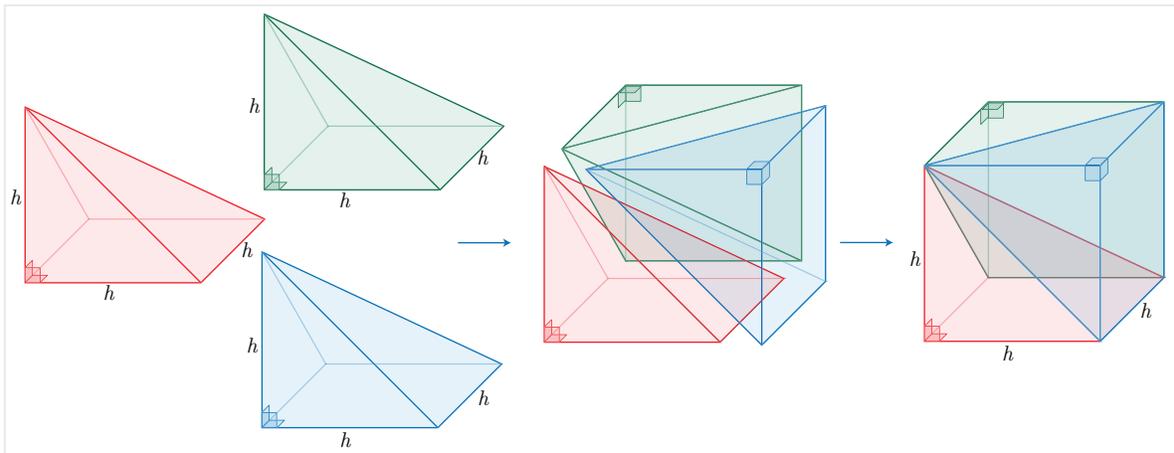}
	\caption{Splitting a cube into three equal pyramids}
	\label{Figure14}
\end{figure}

Consider three copies of a rectangular pyramid of height $h$ whose base is a square of side $h$ too. If we rotate each copy in a certain way and then put them together,   we obtain a cube of dimensions $h,h,h$ as shown in   \cref{Figure14}.  Clearly, the volume of the resulting cube   is equal to the sum of volumes of the  three equal pyramids. Since the volume of the cube  i.e., $h^3$, was known to Elamite scribes, they could  conclude that the volume of a rectangular pyramid with dimensions $h,h,h$ must be   $\frac{1}{3}h^3 $.   

\begin{figure}[H]
	\centering
	\includegraphics[scale=1]{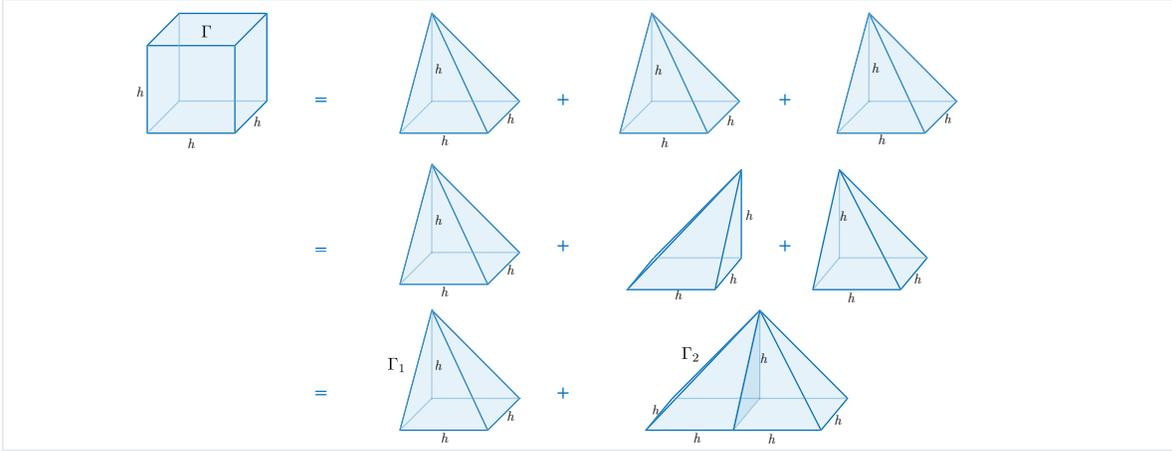}
	\caption{Volume of a rectangular pyramid}
	\label{Figure15}
\end{figure}

 The second fact   now becomes clear.  To see this, take a   cube $\Gamma $ with length $h$, width $h$ and height $h$. As shown   in \cref{Figure14}, this cube  is the union of three copies of      a rectangular pyramid $\Gamma_1$ with equal length, width and height $h$. If we attach two of these copies properly, we get a rectangular  pyramid $\Gamma_2 $ with length $h$, width $2h$ and height $h$ (see \cref{Figure15}).

This means $\Gamma $ is the union of $\Gamma_1$ and $\Gamma_2$ which implies that
\begin{equation}\label{equ-SMT14-ca}
	V_{\Gamma}=V_{\Gamma_1}+V_{\Gamma_2}
\end{equation} 
and also
\begin{equation}\label{equ-SMT14-cb}
	V_{\Gamma}=3V_{\Gamma_1}.
\end{equation} 
Since $ V_{\Gamma}=h^3$, both \cref{equ-SMT14-ca} and \cref{equ-SMT14-cb}  imply that $ V_{\Gamma_2}=h^3-\frac{1}{3}h^3 $ which proves \cref{equ-SMT14-ba}.

Now, we return to the first problem of \textbf{SMT No.\,14}. Lines 1-3 tell us that the volume\index{volume of a grain-heap} of the grain-heap\index{grain-heap},   $V=V_{\Lambda}+V_{P_1}+V_{P_2} $,  is equal to  $14,24$ volume-sar\index{volume-sar (capacity unit)}, so it follows from  \cref{equ-SMT14-ab}   that 
\begin{equation}\label{equ-SMT14-d}
	x h^2  + 2(1-0;20)h^3 = 14,24.
\end{equation}

\noindent
Since\footnote{Note that 1 {\fontfamily{qpl}\selectfont nindan}\index{nindan (length unit)} is equal to 12 k\`{u}\v{s}.}\index{kuz@k\`{u}\v{s} (length unit)}
$$1\ \text{{\fontfamily{qpl}\selectfont k\`{u}\v{s}}}= 0;5\ \text{{\fontfamily{qpl}\selectfont nindan}\index{nindan (length unit)}} $$
and  
\[1\ \text{{\fontfamily{qpl}\selectfont volume-sar}\index{volume-sar (capacity unit)}} = \left( 1\ \text{{\fontfamily{qpl}\selectfont nindan}\index{nindan (length unit)}}^2\right)   \times \left( 1\ \text{{\fontfamily{qpl}\selectfont k\`{u}\v{s}}}\right) = 0;5\  \text{{\fontfamily{qpl}\selectfont nindan}\index{nindan (length unit)}}^3\]
so     according to   lines 3-4,  
\[ 14,24\ \text{{\fontfamily{qpl}\selectfont volume-sar}\index{volume-sar (capacity unit)}} = (14,24)\times(0;5)\ \text{{\fontfamily{qpl}\selectfont nindan}\index{nindan (length unit)}}^3 = 1,12\  \text{{\fontfamily{qpl}\selectfont nindan}\index{nindan (length unit)}}^3.\]
According to lines 5-14, in order to compute the value of $x$, we can substitute\index{substitution method} $h=3$ {\fontfamily{qpl}\selectfont nindan}\index{nindan (length unit)} in  \cref{equ-SMT14-d} and simplify. Note that we have converted all the units involved into {\fontfamily{qpl}\selectfont nindan}. We have 
\begin{align*}
	&~~  x h^2  + 2(1-0;20)h^3 = 14,24 \\
	\Longrightarrow~~&~~    3^2x  + 2(1-0;20)3^3 = 1,12 \\
	\Longrightarrow~~&~~    9x  +  (1;20)\times 27 = 1,12 \\
	\Longrightarrow~~&~~ 9x  +  36 = 1,12 \\  
	\Longrightarrow~~&~~ 9x=  1,12-36\\
	\Longrightarrow~~&~~  9x=   36\\
	\Longrightarrow~~&~~   x=  \frac{1}{9}\times 36\\
	\Longrightarrow~~&~~   x=  (0;6,40)\times 36\\
	\Longrightarrow~~&~~ x= 4. 
\end{align*}   
Finally, by lines 16-17,  we can use    \cref{equ-SMT14-a} to write 
\[  y = 2\times3 = 6\]
and
\[ z = 4+(2\times 3)=4+6  = 10.\]

In the second problem the {\fontfamily{qpl}\selectfont volume-sar}\index{volume-sar (capacity unit)}  14,24 is converted to the {\fontfamily{qpl}\selectfont s\`{i}la}\index{szila@s\`{i}la (capacity unit)}  unit by the storage constant\index{storage constant} 8,0,0\footnote{In fact, 1 {\fontfamily{qpl}\selectfont volume-sar}\index{volume-sar (capacity unit)} is equal to 5,0,0 {\fontfamily{qpl}\selectfont s\`{i}la}\index{szila@s\`{i}la (capacity unit)} (or 18,800 liters). The constant here is a bit larger than this for unknown reasons!}:
\[14,24~\text{{\fontfamily{qpl}\selectfont volume-sar}}~=~(8,0,0)\times(14,24)~\text{{\fontfamily{qpl}\selectfont s\`{i}la}\index{szila@s\`{i}la (capacity unit)}} = 1,55,12,0,0~\text{{\fontfamily{qpl}\selectfont s\`{i}la}\index{szila@s\`{i}la (capacity unit)}}\]
and further to
\[23~\text{{\fontfamily{qpl}\selectfont gur}}_7~\text{and}~2,24~\text{{\fontfamily{qpl}\selectfont gur}\index{gur (capacity unit)}}\]
where   {\fontfamily{qpl}\selectfont gur$_7 $} is the largest capacity unit such that 1 {\fontfamily{qpl}\selectfont gur$_7 $} = 5,0,0,0  {\fontfamily{qpl}\selectfont s\`{i}la}\index{szila@s\`{i}la (capacity unit)} in the Old Babylonian\index{Old Babylonian} period and as we saw before 1 {\fontfamily{qpl}\selectfont gur}\index{gur (capacity unit)} =5,0 {\fontfamily{qpl}\selectfont s\`{i}la}\index{szila@s\`{i}la (capacity unit)}.

\section{Conclusion}
Clearly,  the Elamite scribes like their Babylonian counterparts  were familiar with the basics of solid geometry and    knew how to compute the volume of  three-dimensional figures such as  cubes, prisms and truncated pyramids.   

Our   mathematical interpretation of \textbf{SMT No.\,14}  reveals that  the scribes of Susa  used a formula for the volume of a rectangular pyramid which enables them to calculate  the right value. It confirms that they knew the volume formula for a pyramid  even if they may  not have  expressed it explicitly and this ability on their part    is of considerable interest to those researching   the history of mathematics.

\section*{Appendix: Mathematical Tablet BM 85194}\phantomsection\label{appendix}
In this appendix, we give the transliteration, the translation and the mathematical interpretation of lines 41-49 on the reverse of  \textbf{BM 85194}.

\subsection*{Transliteration}

\begin{note1} 
	\underline{Reverse II:  Lines 41-49}\\
	(L41)\hspace{2mm}  \textit{\textbottomtiebar{h}i-ri-tum} 10-ta-àm \textit{mu-\textbottomtiebar{h}u} 18 sukud \textit{i-na} 1 kùš 1 šà-gal\\
	(L42)\hspace{2mm}  \textit{sà-súm ù} sa\textbottomtiebar{h}ar-\textbottomtiebar{h}i-a za-e 5 \textit{ù} 5 ul-gar 10 \textit{ta-mar} \\
	(L43)\hspace{2mm}  [10] \textit{a-na} 18 sukud \textit{i-ši} 3 \textit{ta-mar} 3 \textit{i-na} 10 ba-zi 7 \\
	(L44)\hspace{2mm}  [\textit{ta-mar}] \textit{sà-súm} nígin-na \textit{sà-súm ù mu-\textbottomtiebar{h}u} 10 ul-gal 17 \textit{ta-mar}\\
	(L45)\hspace{2mm}   [$\frac{1}{2}$ 17 \textit{\textbottomtiebar{h}e-pé}] 8,30 \textit{ta-mar} nigin 1,12,15 \textit{ta-mar}\\
	(L46)\hspace{2mm}  1,12,[15 gar]-ra igi-2-gál 3 dirig \textit{ša mu-\textbottomtiebar{h}u} ugu \\
	(L47)\hspace{2mm} \textit{sà-súm} nig[in šušana] 45 \textit{a-na} 1,12,15 da\textbottomtiebar{h}-\textbottomtiebar{h}a-\textit{ma} \\
	(L48)\hspace{2mm} 1,13 \textit{ta-mar} 1[8] \textit{a-na} 1,13 \textit{i-ši} 22,30 (sic) \textit{ta-mar}\\
	(L49)\hspace{2mm}  2 (èše) 1 (iku) 1 (ubu) gán sa\textbottomtiebar{h}ar-\textbottomtiebar{h}i-a \textit{ki}-$<$\textit{a-am}$ > $ \textit{ne-pé-šum} 
	
\end{note1}

\subsection*{Translation}

\underline{Reverse II:  Lines 41-49}
\begin{tabbing}
	\hspace{18mm} \= \kill 
	(L41)\> \tabfill{A hole dug in the ground. Each of the sides of the upper surface is 10 ({\fontfamily{qpl}\selectfont nindan} $\approx 60m$) in length. The depth is 18 ({\fontfamily{qpl}\selectfont kùš} $\approx 9m$). The inclination of the slope is 1 ({\fontfamily{qpl}\selectfont kùš}) per 1  {\fontfamily{qpl}\selectfont kùš} ($ \approx 50cm$).}\\
	(L42)\> \tabfill{(What are the sides of) the base surface and the volume? (When) you (perform the operation), add 0;5 and 0;5 together, and you see 0;10.} \\ 
	(L43)\> \tabfill{Multiply 0;10 by 18 of the depth, and you see 3. Subtract 3 from 10, and you see 7,}\\ 
	(L44)\> \tabfill{the side of the base surface. On the other hand, add the side of the base and the side of the upper together, and you see 17.}\\
	(L45)\> \tabfill{Halve 17, and you see 8;30. Square (it), and you see 1,12;15.}\\ 
	(L46)\> \tabfill{Put down 1,12;15. $ \frac{1}{2} $ of 3 that is the difference between the upper side and the base side.}\\
	(L47)\> \tabfill{square (it) and (take) $ \frac{1}{3} $ (of the result). Add 0;45 to 1,12;15, and}\\
	(L48)\> \tabfill{you see 1,13. Multiply 18 by 1,13, (and) you see 22,30 (error for 21,54).}\\ 
	(L49)\> \tabfill{2 ({\fontfamily{qpl}\selectfont èše}) 1 ({\fontfamily{qpl}\selectfont iku}) 1 ({\fontfamily{qpl}\selectfont ubu}) ($ = 22,30 $ {\fontfamily{qpl}\selectfont volume-sar}) is the volume (of the hole). Such is the procedure.}
	
\end{tabbing}

\subsection*{Mathematical Interpretation}
\begin{figure}[H]
	\centering
	\includegraphics[scale=1]{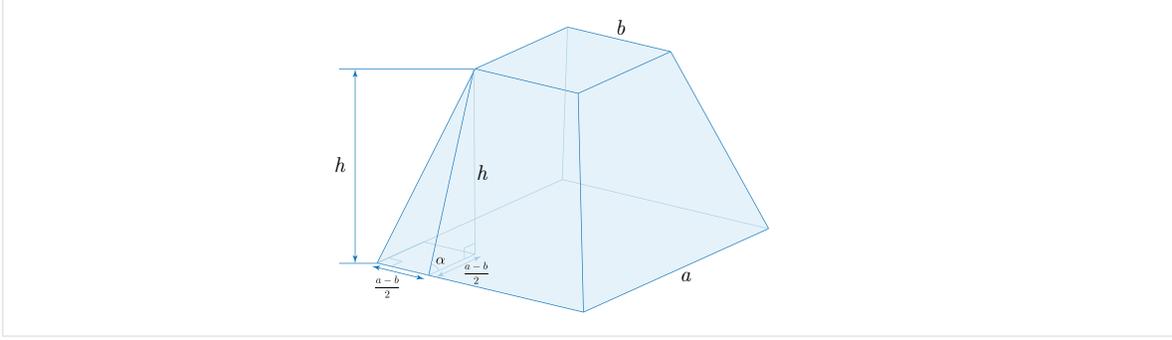}
	\caption{Volume of a pyramidal frustum according to \textbf{BM 85194}}
	\label{Figure9-4}
\end{figure}

Here, the scribe is considering a pyramidal frustum whose base and top are squares with sides $a $ and $b$ and the height (depth) is $h$ (see \cref{Figure9-4}). According to the date in text, $a=10$ {\fontfamily{qpl}\selectfont nindan} and  $h=18$ {\fontfamily{qpl}\selectfont kùš}. Since the slope is 1, the angle $\alpha$ in the figure must be $45^{\circ}$, so 
\[ 1=\tan(45^{\circ})=\tan(\alpha)=\frac{h}{\frac{a-b}{2}} \]
and thus
\[ h=\frac{a-b}{2}.  \]

Since $a=10$ {\fontfamily{qpl}\selectfont nindan},   $h=18$ {\fontfamily{qpl}\selectfont kùš}, and 1 {\fontfamily{qpl}\selectfont kùš} = 0;5 {\fontfamily{qpl}\selectfont nindan}, we get
\[ a-b=2\times (0;5)\times 18=3~\text{{\fontfamily{qpl}\selectfont nindan}}.  \]
This implies that
\[ b=a-(a-b)=10-3=7~\text{{\fontfamily{qpl}\selectfont nindan}} \]
  and 
\[ a+b=7+10=17~\text{{\fontfamily{qpl}\selectfont nindan}}. \]  
  Therefore
\[ \left(\frac{a+b}{2}\right)^2= \left(\frac{17}{2}\right)^2=(8;30)^2=1,12;15~\text{{\fontfamily{qpl}\selectfont nindan}}^2\]  
and
\[ \frac{1}{3}\left(\frac{a-b}{2}\right)^2= \frac{1}{3}\left(\frac{3}{2}\right)^2=\frac{1}{3}\times(1;30)^2=0;45~\text{{\fontfamily{qpl}\selectfont nindan}}^2.\] 
The volume is
\begin{align*}
	V&=\left[\left(\frac{a+b}{2}\right)^2+\frac{1}{3}\left(\frac{a-b}{2}\right)^2\right]h\\
	&=18\times  (1,12;15 + 0;45)\\
	&=18\times (1,13)\\
	&=21,54~\text{{\fontfamily{qpl}\selectfont volume-sar}}. 
\end{align*}

\end{document}